\documentclass[12pt]{article}
\usepackage{amsmath,amssymb,amsfonts,amsthm}
\newenvironment{myabstract}{\par\noindent
{\bf Abstract . } \small } {\par\vskip8pt minus3pt\rm}
\newcounter{item}[section]
\newcounter{kirshr}
\newcounter{kirsha}
\newcounter{kirshb}
\newenvironment{enumarab}{\setcounter{kirshb}{1}
\begin{list}{(\arabic{kirshb})}{\usecounter{kirshb}} }{\end{list}}

\newenvironment{athm}[1]{\vskip3mm\par\noindent
{\bf #1 }. \slshape } {\upshape\par\vskip10pt minus3pt}
\newenvironment{demo}[1]{\noindent{\bf #1.}\upshape\mdseries}
{\nopagebreak{\hfill\rule{2mm}{2mm}\nopagebreak}\par\normalfont}
\theoremstyle{definition}

\def\CA{{\bf CA}}
\def\RCA{{\bf RCA}}

\def\(R)RA{{\bf (R)RA}}
\def\RA{{\bf RA}}
\def\RRA{{\bf RRA}}

\def\Nr{{\mathfrak{Nr}}}

\def\CA{{\bf CA}}

\def\RCA{{\bf RCA}}

\def\(R)RA{{\bf (R)RA}}
\def\RA{{\bf RA}}
\def\RRA{{\bf RRA}}

\def\Nr{{\mathfrak{Nr}}}

\def\CA{{\bf CA}}

\def\RCA{{\bf RCA}}

\def\(R)RA{{\bf (R)RA}}
\def\Ra{{\bf Ra}}
\def\RA{{\bf RA}}
\def\RRA{{\bf RRA}}

\def\Cm{{\mathfrak{Cm}}}
\def\Tm{{\mathfrak{Tm}}}

\title{Omitting types and complete representations}
\author{Tarek Sayed Ahmed}

\title{A remark on completions of the class $S\Nr_n\CA_{n+k},$ $k\geq 3$} 
\author{Tarek Sayed Ahmed\\
Department of Mathematics, Faculty of Science,\\ 
Cairo University, Giza, Egypt.
  }
\begin{document}
\maketitle
\begin{myabstract} We give a sufficient condition that implies that for for $k\geq n+3$, the class 
$S\Nr_n\CA_{n+3}$ is not closed under
completions and is not Sahlqvist axiomatizable. We compare this condition to existing results in the literature.

\end{myabstract}

\bibliographystyle{plain}

We use a construction of Hirsch and Hodkinson \cite{HHbook}. We follow the notation in op.cit. 
$\RA$ stands for the class of relation algebras and
$\CA_n$ stands for the class of cylindric algebras of dimension $n$.
$\Ra\CA_n$ stands for the class of relation algebra reducts of $\CA_n$.
$\RCA_n$ stands for the class of representable $\CA_n$'s and for $k<n$, $\Nr_k\CA_n$ stands for 
the class of neat $k$ reducts of algebras in $\CA_n$.
It is known that $\RCA_n={\bf S}\Nr_n\CA_{\omega}$ for any finite $n$ and that for $k\in \omega$ and $n>2$, 
$\RCA_n\subset {\bf S}\Nr_n\CA_{n+k+1}\subset 
{\bf S}\Nr_n\CA_{n+k}$ \cite{HHbook}. 
Let $\RA_n$ be the class of subalgebras of atomic relation algebras having $n$ dimensional relational basis.
Then ${\bf S}\Ra\CA_n\subseteq \RA_n$ \cite{HHbook}. 
The full complex algebra of an atom structure $S$
will be denoted by $\Cm S$, and the term algebra by $\Tm S.$
$S$ could be a relation atom structure or a cylindric atom structure.

\begin{athm}{Theorem 1} Let $n\geq 3$. Assume that for any simple atomic relation algebra $\cal A$ with atom structure $S$, 
there is a cylindric atom structure $H$ such that:
\begin{enumarab}
\item If $\Tm S\in \RRA$, then $\Tm H\in \RCA_n$.
\item $\Cm S$ is embeddable in $\Ra$ reduct of $\Cm H$.
\end{enumarab}
Then for all $k\geq 3$, $S\Nr_n\CA_{n+k}$ is not closed under completions.

\end{athm}
\begin{demo}{Proof} Let $S$ be a relation atom structure such that $\Tm S$ is representable while $\Cm S\notin \RA_6$.
Such an atom structure exists \cite{HHbook} Lemmas 17.34-17.36 and are finite. 
It follows that $\Cm S\notin {\bf S}\Ra\CA_n$.
Let $H$ be the $\CA_n$ atom structure provided by the hypothesis of the previous theorem.  
Then $\Tm H\in \RCA_n$. We claim that $\Cm H\notin {\bf S}\Nr_n\CA_{n+k}$, $k\geq 3$. 
For assume not, i.e. assume that $\Cm H\in {\bf S}\Nr_n\CA_{n+k}$, $k\geq 3$.
We have $\Cm S$ is embeddable in $\Ra\Cm H.$  But then the latter is in ${\bf S}\Ra\CA_6$
and so is $\Cm S$, which is not the case.
 \end{demo}
 
\begin{athm}{Corollary 2} Assume the hypothesis in Theorem 1. Then the following hold:

\begin{enumarab}

\item There exist two atomic 
cylindric algebras of dimension $n$  with the same atom structure, only one of which is 
representable.

\item For $n\geq 3$ and $k\geq 3$, ${\bf S}\Nr_n\CA_{n+k}$
is not closed under completions and is not atom-canonical. 
In particular, $\RCA_n$ is not atom-canonical.

\item There exists a non-representable $\CA_n$ with a dense representable
subalgebra.

\item For $n\geq 3$ and $k\geq 3$,  ${\bf S}\Nr_n\CA_{n+k}$ 
is not Sahlqvist axiomatizable. In particular, $\RCA_n$ is not Sahlqvist axiomatizable.

\item There exists an atomic representable 
$\CA_n$ with no complete representation.

\end{enumarab}
\end{athm}
\begin{demo}{Proof} \cite{t}
\end{demo}
Monk and Maddux constructs such an $H$ for $n=3$ and Hodkinson constructs  an $H,$ but $H$ does not satisfy 2 \cite{h}.


\begin{thebibliography}{100}
\bibitem{h} Hokinson {\it Constructing cylindric and polyadic algebras from atomic relation algebras}
\bibitem{HHbook} Hirsch R., Hodkinson.I., {\it Relation algebras by games.}
Studies in Logic and the Foundations of Mathematics. Volume 147.
\bibitem{t} Sayed Ahmed {\it The class ..is not closed under completions}

\end{thebibliography}
\end{document}